\documentstyle[twoside,epsfig]{article}
\newcommand{\Z}{{\bf Z}}
\oddsidemargin=\evensidemargin
\catcode`\@=11
\long\def\@makefntext#1{
\protect\noindent \hbox to 3.2pt {\hskip-.9pt  
$^{{\eightrm\@thefnmark}}$\hfil}#1\hfill}		

\def\ps@myheadings{\let\@mkboth\@gobbletwo		
\def\@oddhead{\hbox{}
\rightmark\hfil\eightrm\thepage}   
\def\@oddfoot{}\def\@evenhead{\eightrm\thepage\hfil
\leftmark\hbox{}}\def\@evenfoot{}
\def\sectionmark##1{}\def\subsectionmark##1{}}

\catcode`@=11		      		     
\def\ps@plain{\let\@mkboth\@gobbletwo
     \def\@oddhead{}\def\@oddfoot{\eightrm\hfil\thepage
     \hfil}\def\@evenhead{}\let\@evenfoot\@oddfoot}

\newcounter{sectionc}\newcounter{subsectionc}\newcounter{subsubsectionc}
\renewcommand{\section}[1] {\vspace{12pt}\addtocounter{sectionc}{1} 
\setcounter{subsectionc}{0}\setcounter{subsubsectionc}{0}\noindent 
	{\tenbf\thesectionc. #1}\par\vspace{5pt}}
\renewcommand{\subsection}[1] {\vspace{12pt}\addtocounter{subsectionc}{1} 
	\setcounter{subsubsectionc}{0}\noindent 
	{\bf\thesectionc.\thesubsectionc. 
	{\kern1pt \bfit #1}}\par\vspace{5pt}}
\renewcommand{\subsubsection}[1] {\vspace{12pt}
	\addtocounter{subsubsectionc}{1}
	\noindent
	{\tenrm\thesectionc.\thesubsectionc.\thesubsubsectionc.	{\kern1pt 
	\it #1}}\par\vspace{5pt}}

\newcounter{appendixc}
\newcounter{subappendixc}[appendixc]
\newcounter{subsubappendixc}[subappendixc]

\renewcommand{\appendix}[1] {\vspace{12pt}	
	\refstepcounter{appendixc}		
	\setcounter{figure}{0}
	\setcounter{table}{0}
	\setcounter{lemma}{0}
	\setcounter{theorem}{0}
	\setcounter{corollary}{0}
	\setcounter{definition}{0}
	\setcounter{equation}{0}
	\setcounter{proposition}{0}
	\setcounter{conjecture}{0}
	\setcounter{definition}{0}
	\setcounter{remark}{0}
	\renewcommand{\thefigure}{\Alph{appendixc}.\arabic{figure}}
	\renewcommand{\thetable}{\Alph{appendixc}.\arabic{table}}
	\renewcommand{\theappendixc}{\Alph{appendixc}}
	\renewcommand{\theequation}{\Alph{appendixc}.\arabic{equation}}
	\noindent{\tenbf Appendix \theappendixc #1}\par\vspace{5pt}}

\newcommand{\textlineskip}{\baselineskip=13pt}
\newcommand{\smalllineskip}{\baselineskip=10pt}

\newcommand{\copyrightheading}[1]
	{\vspace*{-2.5cm}\smalllineskip{\flushleft
	{\footnotesize Journal of Knot Theory and Its Ramifications #1}\\
   	{\footnotesize \copyright\kern2pt World Scientific 
         Publishing Company}\\
         }}

\newcommand{\publisher}[2]{{\begin{center}\footnotesize\smalllineskip 
	Received: March 10, 2003 #1\\
	Revised: May 28, 2003 #2
        \end{center}
	}}

\def\abstracts#1#2#3#4{{
	\centering{\begin{minipage}{4.5in}\footnotesize\baselineskip=10pt
	\centerline{ABSTRACT} 
	\parindent=15pt #1\par 
	\parindent=15pt #2\par
	\parindent=15pt #3\par
	\parindent=15pt #4\par
	\end{minipage}}\par}} 

\def\keywords#1{{ 
	\centering{\begin{minipage}{4.5in}\footnotesize\baselineskip=10pt
	{\footnotesize\it Keywords}\/: #1
	\end{minipage}}\par}}

\renewenvironment{thebibliography}[1]
	{\frenchspacing
	 \ninerm\baselineskip=11pt
	 \begin{list}{[\arabic{enumi}]}
	{\usecounter{enumi}\setlength{\parsep}{0pt}
	 \setlength{\leftmargin 13.7pt}{\rightmargin 0pt} 
	 \setlength{\itemsep}{0pt} \settowidth
	{\labelwidth}{[#1]}\sloppy}}{\end{list}}

\newcounter{itemlistc}
\newcounter{romanlistc}
\newcounter{alphlistc}
\newcounter{arabiclistc}

\newcommand{\fcaption}[1]{
        \refstepcounter{figure}
        \setbox\@tempboxa = \hbox{\footnotesize Fig.~\thefigure. #1}
        \ifdim \wd\@tempboxa > 5in
           {\begin{center}
        \parbox{5in}{\footnotesize\smalllineskip Fig.~\thefigure. #1}
            \end{center}}
        \else
             {\begin{center}
             {\footnotesize Fig.~\thefigure. #1}
              \end{center}}
        \fi}

\newcommand{\tcaption}[1]{
        \refstepcounter{table}
        \setbox\@tempboxa = \hbox{\footnotesize Table~\thetable. #1}
        \ifdim \wd\@tempboxa > 5in
           {\begin{center}
        \parbox{5in}{\footnotesize\smalllineskip Table~\thetable. #1}
            \end{center}}
        \else
             {\begin{center}
             {\footnotesize Table~\thetable. #1}
              \end{center}}
        \fi}

\def\pmb#1{\setbox0=\hbox{#1}
	\kern-.025em\copy0\kern-\wd0
	\kern.05em\copy0\kern-\wd0
	\kern-.025em\raise.0433em\box0}

\def\fnt#1#2{\footnotetext{\kern-.3em
	{$^{\mbox{\scriptsize #1}}$}{#2}}}

\def\fpage#1{\begingroup
\voffset=.3in
\thispagestyle{empty}\begin{table}[b]\centerline{\footnotesize #1}
	\end{table}\endgroup}

\def\runninghead#1#2{\pagestyle{myheadings}
\markboth{{\protect\footnotesize\it{\quad #1}}\hfill}
{\hfill{\protect\footnotesize\it{#2\quad}}}}

\font\tenrm=cmr10
 
\font\tenbf=cmbx10
\font\bfit=cmbxti10 at 10pt
\font\ninerm=cmr9
\font\nineit=cmti9
\font\ninebf=cmbx9
\font\eightrm=cmr8

\def\@begintheorem#1#2{\trivlist	
	\item[\hskip\labelsep{\bf #1\ #2.}]} 
\def\@opargbegintheorem#1#2#3{\trivlist
	\item[\hskip\labelsep{\bf #1\ #2\ (#3).}]}

\newenvironment{proof}{\begin{trivlist}
	\item[\noindent]{\it Proof.}}{\quad $\qed$\end{trivlist}} 

    	{\setcounter{itemlistc}{0}		
	 \begin{list}{$\bullet$}		
	{\usecounter{itemlistc}			
	 \leftmargin10pt	       
	 \setlength{\parsep}{0pt}
	 \setlength{\itemsep}{0pt}     
	}}{\end{list}}

	{\setcounter{romanlistc}{0}		
	 \begin{list}{$($\roman{romanlistc}$)$}	
	{\usecounter{romanlistc}		
	 \leftmargin18pt 
	 \setlength{\parsep}{0pt}
	 \setlength{\itemsep}{0pt}	
	 \settowidth{\labelwidth}{#1}                          
	}}{\end{list}}

	{\setcounter{enumii}{0}			
	 \begin{list}{$($\alph{enumii}$)$}	
	{\usecounter{enumii}			
	 \leftmargin18pt		
	 \setlength{\parsep}{0pt}
	 \setlength{\itemsep}{0pt}	
	 \settowidth{\labelwidth}{#1}                          
	}}{\end{list}}

\def\qed{\hbox{${\vcenter{\vbox{			
   \hrule height 0.4pt\hbox{\vrule width 0.4pt height 6pt
   \kern5pt\vrule width 0.4pt}\hrule height 0.4pt}}}$}}

\def\theequation{\thesectionc.\arabic{equation}}  
\begin{document}
\setlength{\textheight}{7.7truein}  

\runninghead{Kauffman-Harary conjecture holds for Montesinos Knots}
{Kauffman-Harary conjecture holds for Montesinos Knots}

\normalsize\textlineskip
\thispagestyle{empty}
\setcounter{page}{1}

\copyrightheading{}		    

\vspace*{0.88truein}

\fpage{1}
\centerline{\bf KAUFFMAN-HARARY CONJECTURE HOLDS FOR MONTESINOS KNOTS}
\vspace*{0.37truein}
\centerline{\footnotesize MARTA M. ASAEDA}
\baselineskip=12pt
\centerline{\footnotesize\it Dept. of Mathematics, Univ. of Maryland}
\baselineskip=10pt
\centerline{\footnotesize\it College Park, MD 20742}
\centerline{\footnotesize\it marta@math.umd.edu}

\vspace*{10pt}
\centerline{\footnotesize  J\'OZEF H. PRZYTYCKI}
\baselineskip=12pt
\centerline{\footnotesize\it Dept. of Mathematics, The George Washington 
University}
\baselineskip=10pt
\centerline{\footnotesize\it 2201 G St. NW, 
Washington, DC 20052}
\centerline{\footnotesize\it przytyck@gwu.edu}

\vspace*{10pt}
\centerline{\footnotesize  ADAM S. SIKORA}
\baselineskip=12pt
\centerline{\footnotesize\it Dept. of Mathematics, SUNY at Buffalo}
\baselineskip=10pt
\centerline{\footnotesize\it Buffalo, NY 14260}
\baselineskip=12pt
\centerline{\footnotesize\it and}
\baselineskip=12pt
\centerline{\footnotesize\it Institute for Advanced Study, 
School of Mathematics}
\baselineskip=10pt
\centerline{\footnotesize\it 1 Einstein Dr., Princeton, NJ 08540}
\centerline{\footnotesize\it asikora@buffalo.edu} 

\vspace*{0.225truein}
\publisher{}

\vspace*{0.21truein} 
\abstracts{The Kauffman-Harary conjecture states that for any reduced 
alternating diagram $K$ of a knot with a prime determinant $p,$ every 
non-trivial Fox 
$p$-coloring of $K$ assigns different colors to its arcs. We generalize 
this conjecture by stating it in terms of homology of 
the double cover of $S^3$ branched along a link.
 In this way we extend the scope of the conjecture 
to all prime alternating 
links of arbitrary determinants. We first prove the Kauffman-Harary 
conjecture for pretzel knots and then we generalize our argument to show
the generalized Kauffman-Harary conjecture for all Montesinos links.
Finally, we speculate on the relation
between the conjecture and Menasco's work on incompressible surfaces 
in exteriors of alternating links.}{}{}{}


\vspace*{10pt}

\vspace*{1pt}\textlineskip	
\section{Introduction}          
\vspace*{-0.5pt}

In this paper we consider the following conjecture by
Kauffman and Harary \cite{H-K}.\\

\noindent
{\bf Conjecture 1 (Kauffman-Harary Conjecture)}
{\it Let $\cal D$ be an alternating knot diagram with no nugatory crossings.
If the determinant of $\cal D$ is a prime number $p$ then
every non-trivial Fox $p$-coloring of $\cal D$ assigns different colors to 
different arcs of $\cal D.$}\ \\

In the first section of the paper we prove the Kauffman-Harary conjecture 
for pretzel knots. In the second section we 
generalize the conjecture in terms of homology of 
 the double branched covers of $S^3$ branched along links
and illustrate it  
by examples. In the third section we prove Kauffman-Harary conjecture and
its generalization for Montesinos links. In the last section
we speculate about an approach to the generalized conjecture by 
relating it to Menasco's results on incompressible surfaces in the 
exteriors of alternating links.

\section{Pretzel Knots and Fox coloring}

In this section we prove the Kauffman-Harary conjecture for pretzel 
knots.  We deal with this special case in order to prepare
a more general setting, in which we replace Fox coloring by
homology of the double branched cover of $S^3$ branched along a link.

\noindent
{\bf Definition 2}{\it
\begin{enumerate}
\item [(i)]
We say that a link (or a tangle) diagram is k-colored if every
arc is colored by one of the numbers $0,1,...,k-1$ (forming the
cyclic group $\Z_k$) in such a way that
at each crossing the sum of the colors of the undercrossings is equal
to twice the color of the overcrossing mod $k$; see Fig.1.1.
\item [(ii)]
The set of $k$-colorings of a diagram $\cal D$ forms an abelian group, 
denoted by $Col_k(\cal{D})$.
\end{enumerate} }
\vspace*{.1in}

\centerline{\psfig{figure=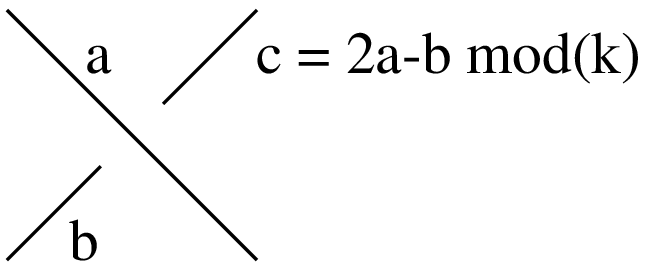,height=1.8cm}}
\centerline{Fig. 1.1}
\vspace*{.2in}

\noindent
{\bf Proposition 3 (Fox)}
{\it
$Col_k({\cal D}) = H_1(M_{\cal{D}}^{(2)},\Z_k) \oplus \Z_k,$ where 
$M_{\cal{D}}^{(2)}$ denotes
the double cover of $S^3$ branched along $\cal D$.}

The first class of knots for which the  Kauffman-Harary conjecture
has been proved is the family of rational (or 2-bridge) knots 
\cite{K-L,PDDGS}. L.~Kauffman challenged us at AMS annual meeting at 
Baltimore in January 2003 to prove the conjecture for pretzel knots and he 
gave some ideas why it should hold \cite{Ka}. 
Initially we 
were skeptical but after analyzing several examples (e.g. pretzel
knot $P(11,7,5,2)$ colored in Fig.1.2) we became convinced that 
the conjecture holds for all alternating knots.

\centerline{\psfig{figure=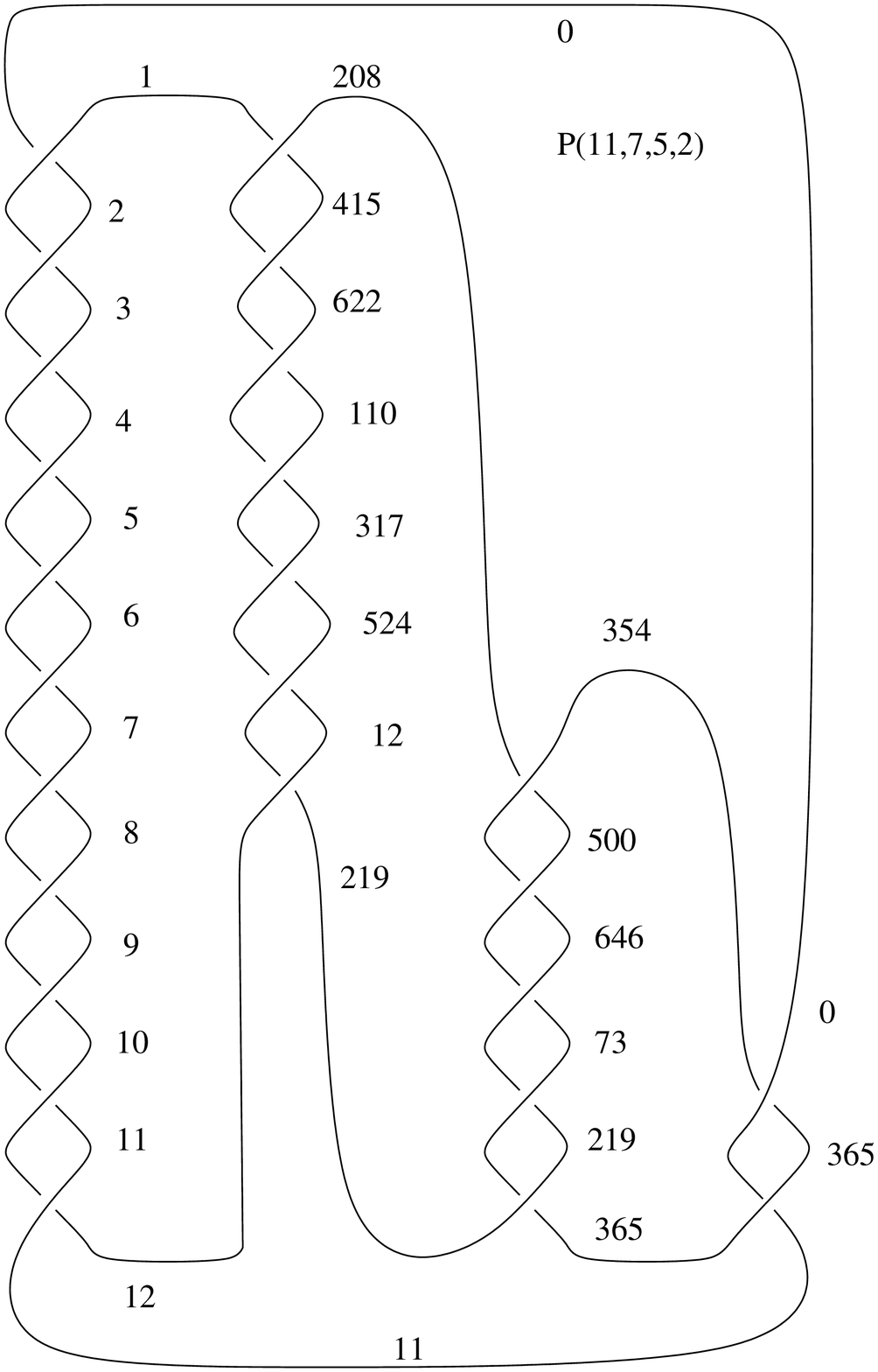,height=10cm}}
\begin{center}
Fig. 1.2. Pretzel knot $K=P(11,7,5,2)$ and its Fox $719$-coloring 
($H_1(M^{(2)}_K)=\Z_{719}$).
\end{center}

First we prove the theorem for pretzel knots, Fig.1.3.\\

\centerline{\psfig{figure= 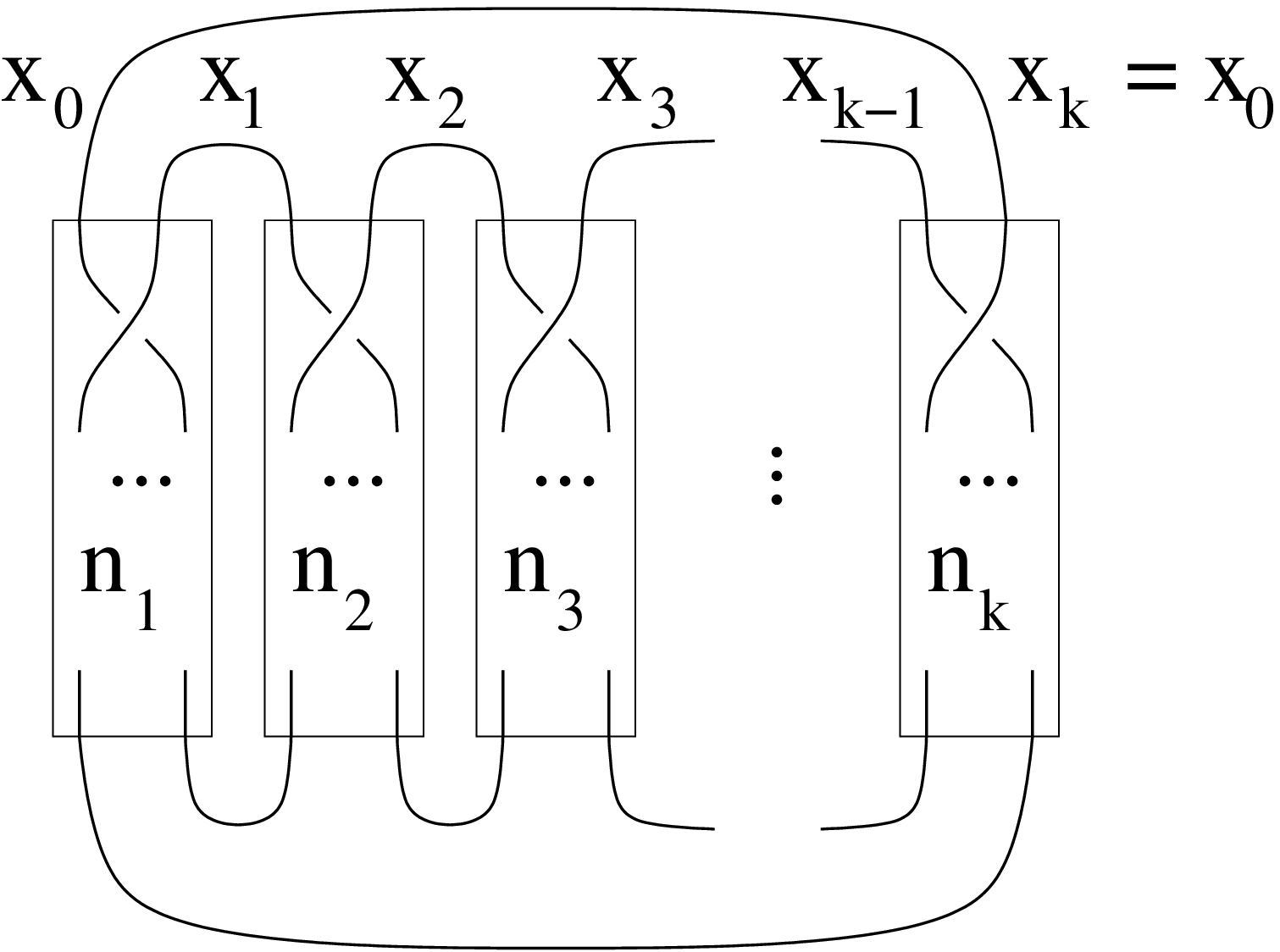,height=5cm}}
\centerline{$n_1z_1=n_2z_2=...=n_kz_k$ mod $D$;\ \ $z_i=x_i-x_{i-1}$.}
\begin{center}
Fig. 1.3. Pretzel knot $ K=P(n_1,n_2,...,n_k)$
\end{center}

\noindent
{\bf Theorem 4}
{\it The Kauffman-Harary conjecture holds for any alternating pretzel knot
diagram, $P(n_1,n_2,...,n_k)$.}

\begin{proof}
Without loss of generality we can assume that $n_1,n_2,...,n_k >0$;
the case $n_1,n_2,...,n_k < 0$ follows in a similar manner.
The determinant of $P(n_1,n_2,...,n_k)$ is $D= 
\Sigma_{i=1}^k n_1n_2...n_{i-1}n_{i+1}...n_k$, by Conway's formula 
\cite{Co}. Let us assume that $D$ is a prime number.
Denote colors of the maximal arcs by $x_0,x_1,...,x_{k-1},
x_k=x_0$, as shown in Fig.1.3. Let $f$ be a nontrivial $D$-coloring. 
We will show that no distinct arcs in two different columns, say $i$th 
and $j$th, use the same color. 
We can assume that $i=1$ and $1<j <k$. Since $D$ is a prime number,
hence any two nontrivial $D$-colorings $f_1,f_2$ are linearly related \ 
$f_2 = a_0 + a_1f_1$, where $a_0$ and $a_1$ are constant (trivial)
colorings. Thus we only analyze a coloring with $x_0=0$ and $x_1=1$,
without loss of generality.
Let $z_i = x_i - x_{i-1}$ (so $z_1=1$). Comparing colorings of
minima of neighboring columns we obtain:
$$n_1z_1 \equiv n_2z_2 \equiv n_3z_3 \equiv ...\equiv n_jz_j
\equiv ...\equiv n_kz_k\ \  {\rm mod\ } D.$$
These equalities determine the color of each $x_i$ uniquely,  
given $x_0 =0$ and $x_1=1$.
With our assumption the first column uses the colors $0,1,2,...,n_1,n_1+1$,
and the $j$th column uses the colors $x_{j-1}$, $x_j=x_{j-1}+z_j$,
$x_{j-1}+2z_j$,...,$x_{j-1} + n_jz_j$, $x_{j-1} + (n_{j}+1)z_j$.
Suppose that an arc in the first column and an arc in the $j$th
column have the same color. Then we have  
$a \equiv x_{j-1} + bz_j$ mod $D$, $0 \leq a \leq n_1 +1$, 
$0 \leq b \leq n_j$.
Observe that
$x_{j-1} = z_{j-1} + z_{j-2}... +
z_2 +1$.
Multiplying both sides of the expression $a \equiv x_{j-1} + bz_j$ mod $D$ by
$n_{j}n_{j-1}...n_2$ one gets
$$n_{j}n_{j-1}...n_2a \equiv n_{j-1}...n_2n_1 b + n_jn_{j-2}....n_2n_1 +...+
n_jn_{j-1}...n_2.$$ Therefore
$$n_{j}n_{j-1}...n_2(a-1) \equiv n_{j-1}...n_2n_1 b + 
n_jn_{j-2}....n_2n_1+...+ n_jn_{j-1}...n_3n_1.$$

Since both sides of the equality are smaller than $D$, it should
hold in $\Z$.  Let us note that  $gcd(n_i,n_j)=1$ for $i \neq j$ since 
$D$ is prime. Note that $n_2$ must divide the right hand sides 
of this equality, which is impossible unless $j=2$.
For $j=2$ we have $(a-1)n_2 \equiv bn_1$, so
$(a-1)$ is divisible by $n_1$ and $b$ is
divisible by $n_2$. This may hold only if $a=n_1+1$ and $b=n_2$ and,  
in this case, the two arcs coincide as minima of neighboring columns.
This completes the proof of Theorem 4. 
\end{proof}

\section{The Generalized Kauffman-Harary Conjecture}

It was noticed in \cite{K-L} that Kauffman-Harary conjecture holds
for any rational (2-bridge) knot without restrictions on the 
determinant of the knot.
However, the formulation of the conjecture needs to be changed in this setting 
from ``every nontrivial $D$-coloring..." to 
``there exists a $D$-coloring...". 

The coloring of $2$-bridge tangle of type $\frac{m}{n}  (m > n)$ with
the maxima colored by $0$ and $1$ is illustrated in Fig.2.1. Observe
that the color of each arc strictly increases as one goes down 
along the diagram. As we close the tangle without introducing any 
new crossings (and not creating a nugatory crossing), we obtain
 a $2$-bridge link of type $\frac{m}{n}$ with determinant $m$ (by Conway's
formula) where $ m > n$. From Fig.2.1 we see that the colors of
the arcs do not exceed $m$ except for the last minimum $m+n \equiv n$
mod $m$, thus the conjecture holds for $2$-bridge links.
Note that in this case
$H_1(M^{(2)}_{L_{\frac{m}{n}}}) =\Z_m$, so the first homology group is cyclic.
\\
\ \\
\centerline{\psfig{figure= 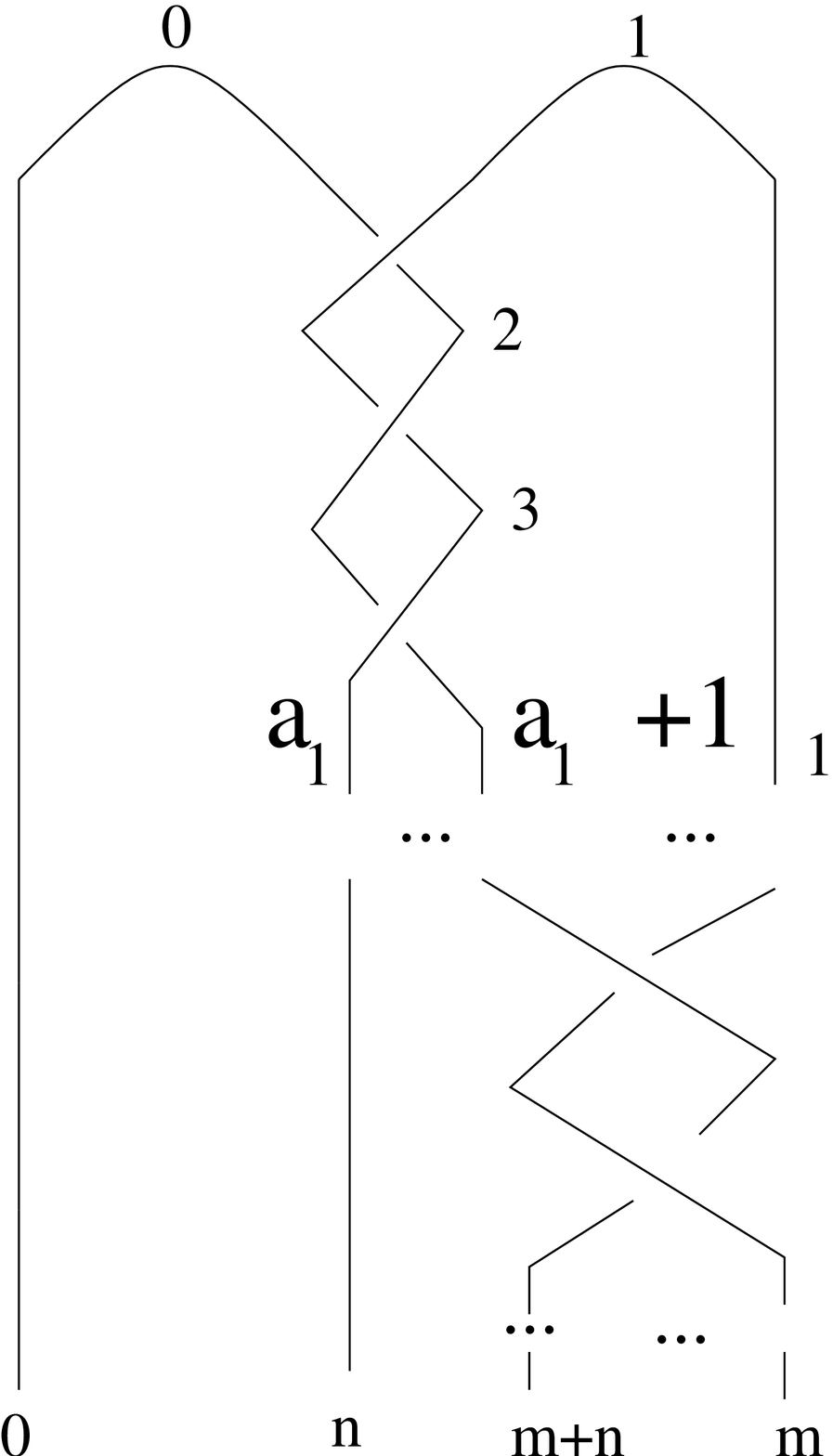,height=8cm}}
\begin{center} 
$\frac{m}{n} = a_k + \frac{1}{a_{k-1} + ... + \frac{1}{a_1}}$ 
\end{center}
\begin{center}
Fig. 2.1. Rational $\frac{m}{n}$-tangle.
\end{center}

We will later use the consequence
of this ``propagation down" proof. Notice that if we color maxima by 
$x_0$ and $x_1$ in place of $0$ and $1$, then an arc colored before by
$c$ is now colored by $x_0+x_1c$.  In Section 3 we consider rational tangles 
of type $\frac{m_i}{n_i}$, where  $m_i \leq n_i $, in which case, if 
$$\frac{m_i}{n_i} = a_{k_i,i} + \frac{1}{a_{k_i-1,i} +...+ 
\frac{1}{a_{1,i}}}$$
then $a_{k_i,i}=0$. Our convention for the diagram of the rational 
$\frac{m}{n}$-tangle (after Conway \cite{Co}) is
presented in Fig.2.2.

\ \\
\centerline{\psfig{figure=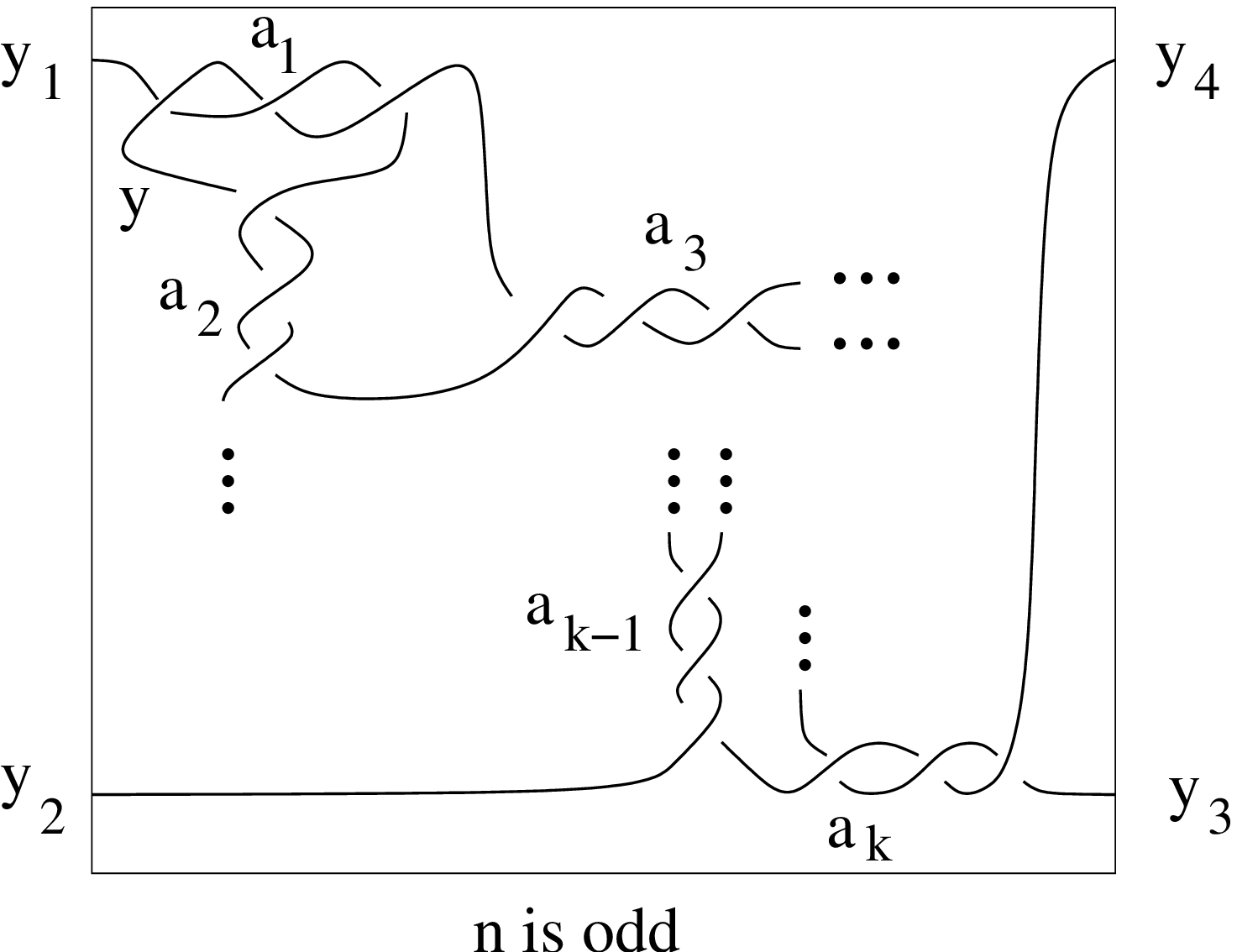,height=5.3cm}\ \ \
\psfig{figure=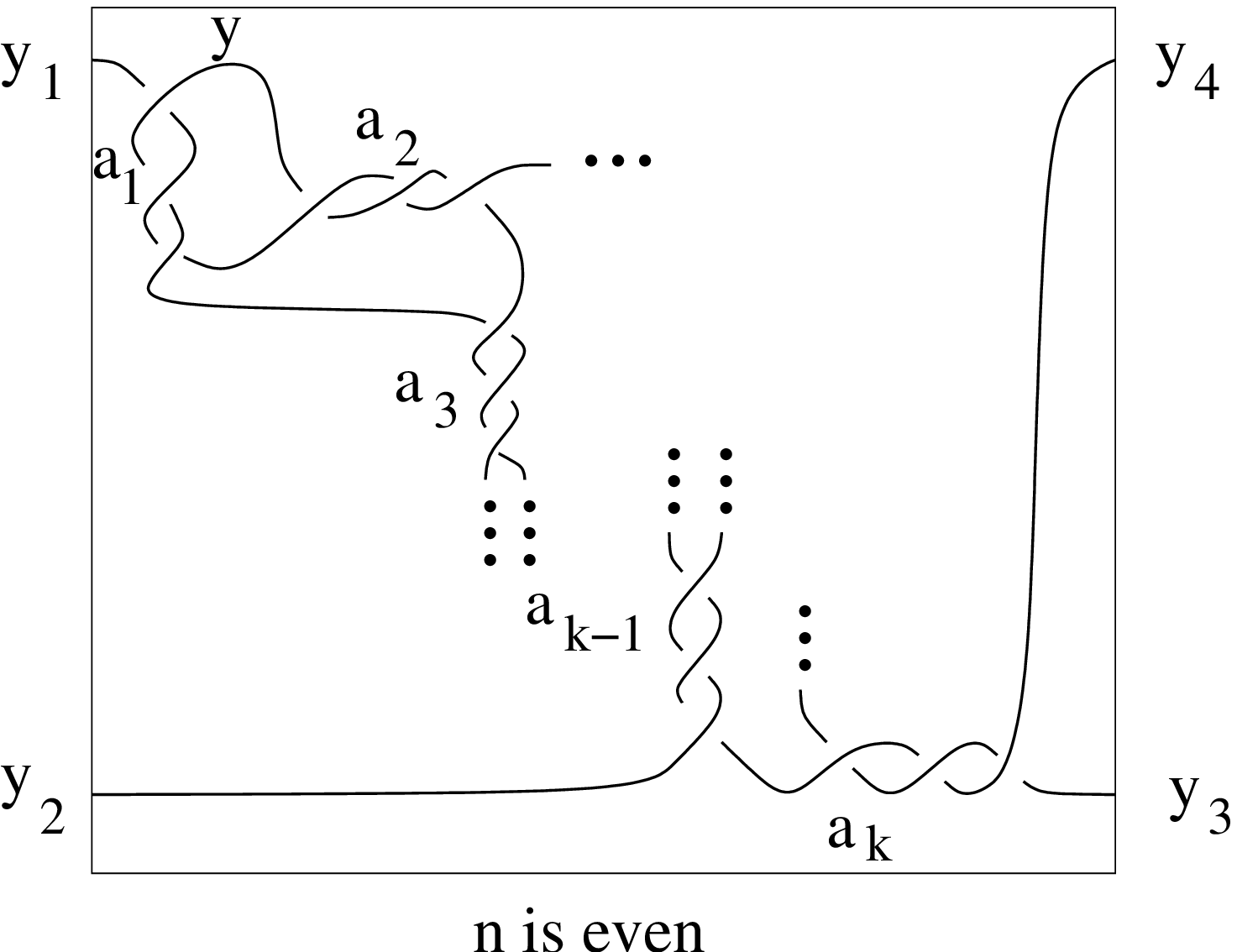,height=5.3cm}}
\begin{center}
The ``internal" maximum is colored by $y$ and then $y_2= m(y-y_1) + y_1$, 
$y_3= (m+n)(y-y_1) + y_1$, and $y_4= n(y-y_1) + y_1$.
\end{center}
\centerline{Fig. 2.2}
\ \\

The pretzel knot $P(15,10,6)$ has cyclic homology 
$H_1(M_{P(15, 10, 6)}^{(2)})= \Z_{300}$. The coloring of this knot
using different colors for each arc is illustrated in Fig.2.3.\\

\centerline{\psfig{figure= 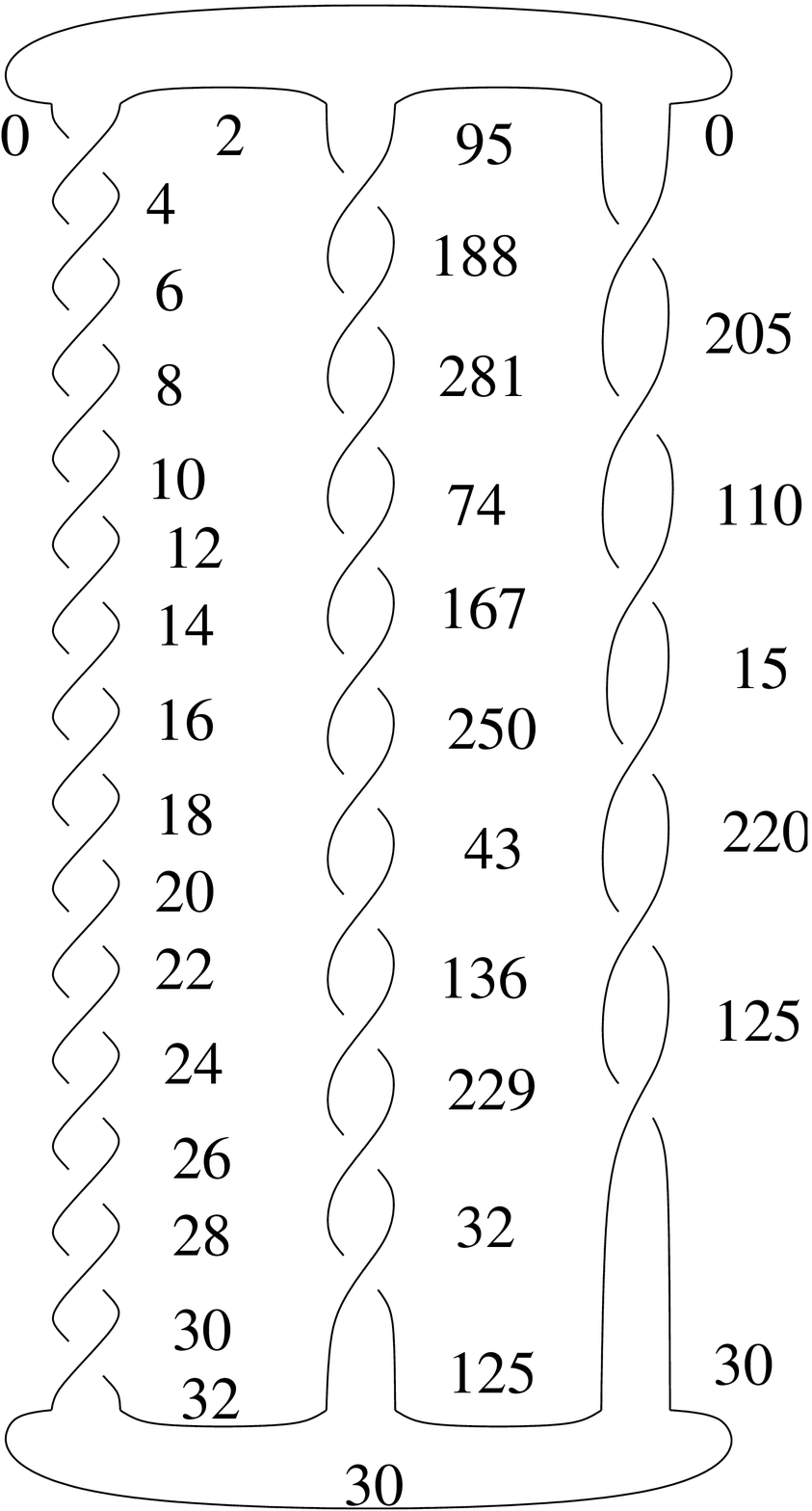,height=9cm}}
\begin{center}
Fig. 2.3. 300-coloring of the pretzel knot, $P(15,10,6)$.
\end{center}

On the other hand the pretzel link $P(3,3,3)$ has the determinant $D=27$ and
does not allow a $27$-coloring with every arc using a different color.
Note that in this case the group $H_1(M_{P(3,3,3)}^{(2)})=\Z_9 \oplus \Z_3$ 
is not cyclic.

While the Kauffman-Harary conjecture requires the determinant of the 
knot to be a prime number, the above examples suggest that a weaker 
requirement could 
suffice. Namely, one could merely assume the homology group 
$H_1(M_{L}^{(2)},\Z)$ to be cyclic (equal to $\Z_D$, where $D$ is not 
necessary a prime number). The conjecture can be
further extended by allowing the elements of the homology group to serve 
as colors. We checked the link $P(3,3,3)$ which has 
$H_1(M_{P(3,3,3)}^{(2)})=\Z_9 \oplus \Z_3$ and the extended conjecture 
holds for this example. We give more details below.

If we decorate arcs of the diagram $L$ by commutative
variables and we quotient the resulting free abelian group generated
by these variables by relations of type $2a-b-c$ 
for every crossing, we get 
$Col(L)= H_1(M_L^{(2)},\Z) \oplus \Z$.
If we choose one arc to be decorated by $0$ then 
we obtain the homology group 
$H_1(M_L^{(2)},\Z)$.
The group $Col_n(L)$ of Fox n-colorings is the module which is 
$\Z_n$-dual to $Col(L)$
(or, equivalently it is the cohomology group with one additional 
$\Z_n$ factor, $H^1(M_L^{(2)},\Z_n) \oplus \Z_n$).
Motivated by this, we suggest that the proper conjecture, generalizing the 
Kauffman-Harary conjecture,
should be the following.
\\ \\
\noindent
{\bf Conjecture 5 (The Generalized Kauffman-Harary (GKH) Conjecture)}\ \\
{\it If $L$ is an alternating diagram of a prime
link without nugatory crossings then different arcs of $L$ represent
different elements of $H_1(M_L^{(2)},\Z)$.}

\noindent
{\bf Remarks 6}\ \\
(1) For a knot $K$ with $H_1(M_K^{(2)},\Z)=\Z_p$, the GKH conjecture 
is equivalent to the Kauffman-Harary conjecture.\\
(2) For a pretzel link $L=P(n_1,n_2,...,n_k)$,  
    the group $H_1(M_L^{(2)},\Z)$ is cyclic if and only if 
    $gcd\{n_1\cdots n_{i-1}n_{i+1}\cdots n_{j-1}n_{j+1}\cdots n_k\ | \ 
    1\leq i < j \leq k \}=1$; 
    compare Proposition 7.\\
(3) Fig.2.4 illustrates the fact that the GKH conjecture holds 
    for the pretzel link $P(3,3,3)$. 
    $H_1(M_{P(3,3,3)}^{(2)},\Z)=\Z_9 \oplus \Z_3$,
    and if we color maxima by $0, x_1, x_2$ respectively, we 
    obtain $x_1$ as a
    generator of $\Z_9$ part and $x_1+x_2$ as a generator of $\Z_3$ part.
    Let us denote $x_1+x_2$ by $u$.
    Then our 9 arcs use 9 different ``colors":
    $0, x_1, 2x_1, 3x_1, 4x_1, 8x_1+u, 6x_1 +2u, 2x_1+u, x_1+2u $.\\
(4) The conjecture obviously fails for non-prime knots
    because for connected sums of knots the connecting arcs always represent
     the same element of the homology group.\\
(5) The conclusion of the conjecture is equivalent to the statement that for
    any pair of arcs
    of the diagram there is a Fox coloring distinguishing them.

\ \\
\centerline{\psfig{figure=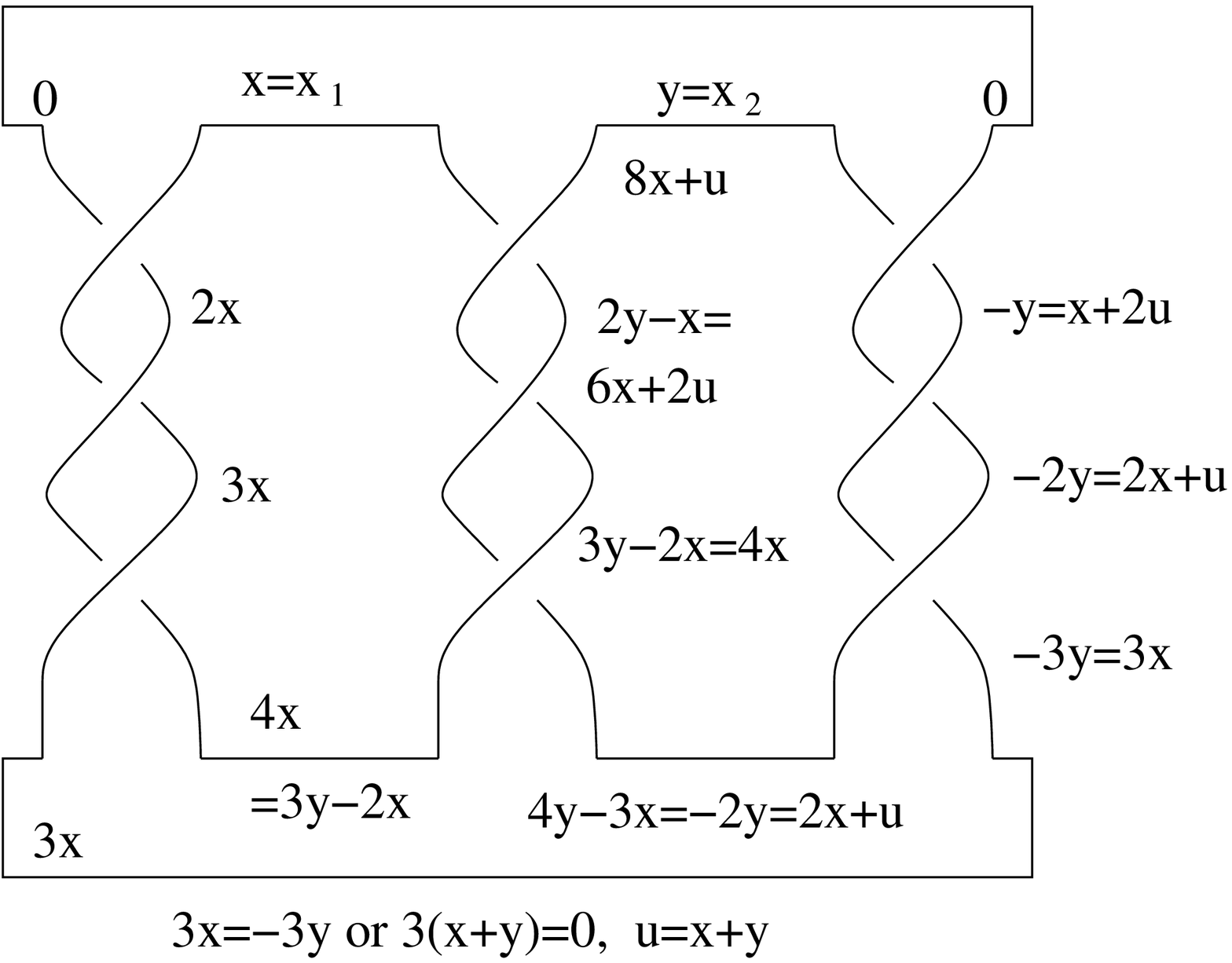,height=7cm}}
\begin{center}
Fig. 2.4. Homology coloring of $P(3,3,3)$.
\end{center}

The first step to prove the GKH conjecture is
to understand the homology group of the double branched cover and
relate it to arc presentation (as in $Col(L)$ group). For the general pretzel
link we have\\

{\bf Proposition 7}
{\it
For the pretzel link $L=P(n_1,n_2,...,n_k)$ the first homology group of
the double branched cover of $S^3$ branched along $L$ has 
the following canonical cyclic decomposition,
$H_1(M_L^{(2)})= \Z_{D_0/D_1} \oplus \Z_{D_1/D_2} \oplus ... 
\oplus \Z_{D_{k-2}},$ where $D_0 = D = 
\Sigma_{i=1}^k n_1\cdots n_{i-1}n_{i+1}\cdots n_k$
$D_1= gcd\{n_1\cdots n_{i-1}n_{i+1}\cdots n_{j-1}n_{j+1}\cdots n_k\}$,..., 
$D_s= gcd\{products \ \ of\  k-s-1\ \  terms\},...,
D_{k-2}=gcd\{n_1,n_2,...,n_k\}$. }

The proof in the more general setting of Montesinos links is given in 
Section 3 (Proposition 8).

\section{Montesinos links}

In this section we prove the GKH
conjecture for alternating Montesinos links (including pretzel knots). 
We draw their diagrams in
the manner similar to pretzel knots (rational $\frac{m_i}{n_i}$-tangle
in place of a column which can be thought as $\frac{1}{n_i}$
rational tangle). Since we deal with alternating Montesinos links, we can
assume that $0< m_i \leq n_i$ and $gcd(m_i,n_i)=1$; see Fig. 3.1 and 3.2.
\\
\ \\
\centerline{\psfig{figure= 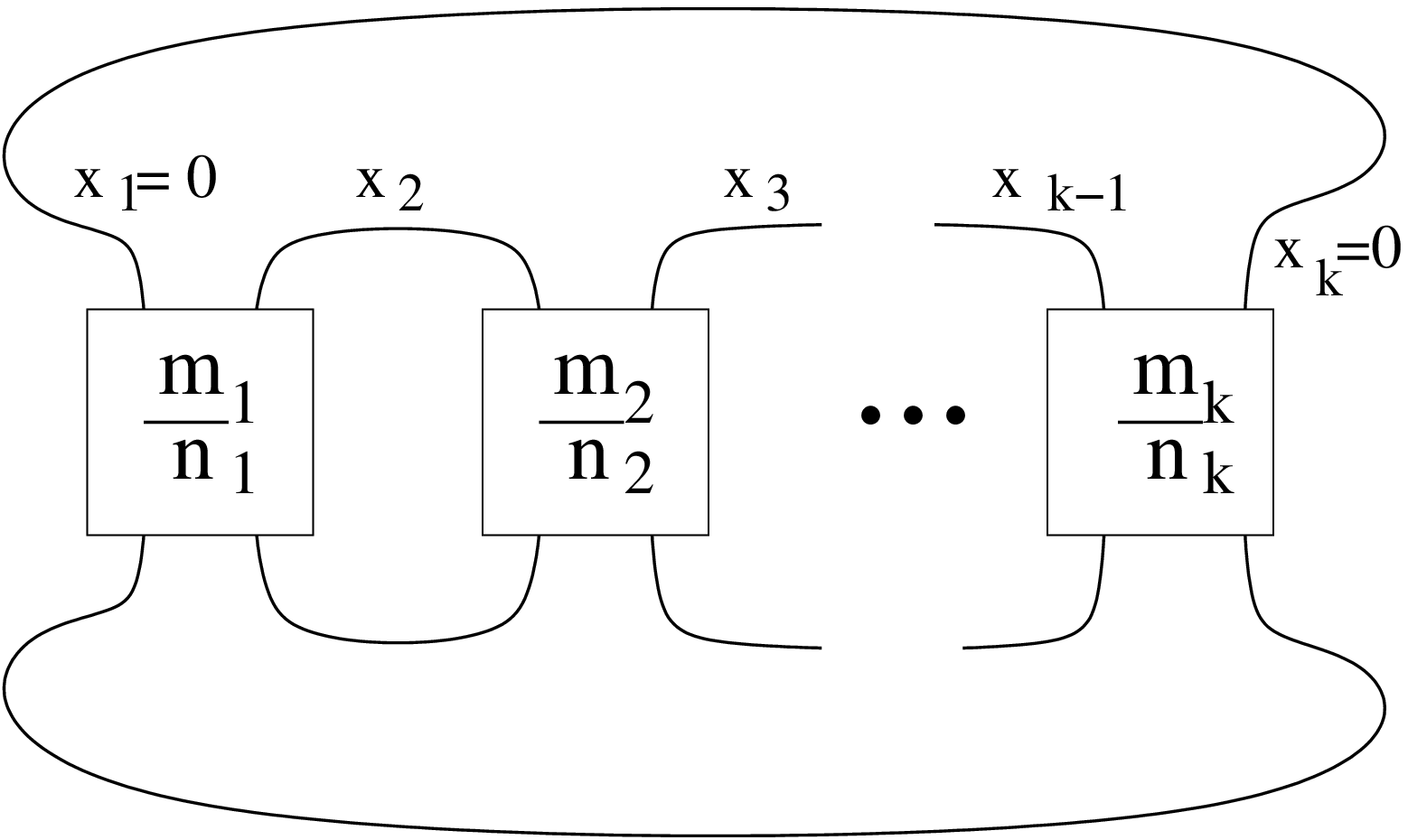,height=5cm}}
\begin{center}
Fig. 3.1; Alternating Montesinos link, $M(\frac{m_1}{n_1},...,\frac{m_k}{n_k})$.\end{center}
\ \\
\ \\
\centerline{\psfig{figure= 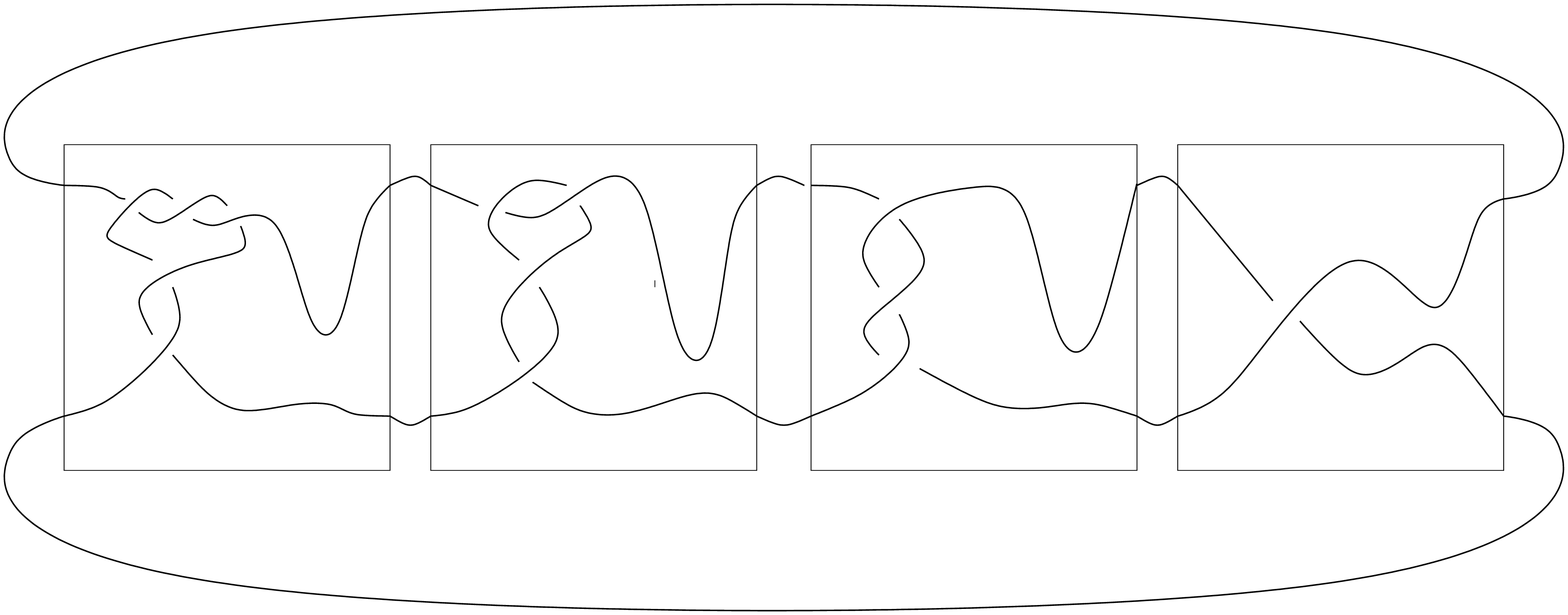,height=4cm}}
\begin{center}
Fig. 3.2;  $M(\frac{3}{7},\frac{2}{5}, \frac{1}{3},\frac{1}{1})$.
\end{center}

Generalizing Proposition 7, we first compute homology
of the double branched cover along a  Montesinos link. Conway's formula
gives the determinant as
$D = \Sigma_{i=1}^k n_1n_2...n_{i-1}m_in_{i+1}...n_k$.\\

{\bf Proposition 8}
{\it The first homology group of the double cover of $S^3$ branched along 
the Montesinos link $L=M(\frac{m_1}{n_1},...,\frac{m_k}{n_k})$
has the following canonical decomposition into cyclic groups,
$H_1(M_L^{(2)})= \Z_{D_0/D_1} \oplus \Z_{D_1/D_2} \oplus ... 
\oplus \Z_{D_{k-2}},$ where $D_0 = D = 
\Sigma_{i=1}^k n_1...n_{i-1}m_in_{i+1}...n_k,$ 
$D_s= gcd (W_{i,s} ),$ and $W_{i,s}$ is the set of products obtained from
$n_1...n_{i-1}m_in_{i+1}...n_k$ by dropping $s$ letters from it. Finally,
$D_{k-1}=1$.}

\begin{proof}
The group $H_{1}(M_L^{(2)},\Z)$ is an abelian group generated by elements
$z_1,z_2,...z_k$ with relations $$ n_{1}z_{1}=n_{2}z_{2}= \ldots =n_{k}z_{k},
\; m_1z_{1}+ m_2z_{2}+ \ldots +m_kz_{k}=0 .$$  Here $z_i=y_i-x_{i-1},$ where
$x_i$ are ``connecting maxima" in the diagram, Fig.3.1 and $y_i$ are
``internal" maxima of rational tangles (compare Fig.2.1,2.2 and 3.2). We have
$x_i-x_{i-1} = m_i(y_i-x_{i-1}) = m_iz_i.$ The relations $n_1z_1=n_jz_j$ are
obtained by comparing labels of the minima of the diagram.
This presentation is described by the matrix $A_k$, where rows 
represent relations of the group. 
That is, we have $H_{1}(M_L^{(2)},\Z)= \Z^k/{\rm Im}(A_k)$, where 
$A_k: \Z^k \to \Z^k$ is the linear map given by $v \mapsto (v)A_k$.
\newfont{\bg}{cmr10 scaled\magstep4}
\newcommand{\bigzerou}{\smash{\lower1.7ex\hbox{\bg 0}}}

$$A_k= \left[
\begin{array}{cccccc}
n_1 & -n_2 & && & \bigzerou \\
n_1 & 0 & -n_3  &&& \\
n_1 & 0 & 0 & -n_4 && \\
\vdots &  \vdots &  & \ddots & \ddots  & \\
n_1 & 0 & \cdots && 0 &  -n_k \\
m_1 & m_2 & \cdots && m_{k-1} & m_k
\end{array} \right].  $$
The canonical decomposition of the group into cyclic groups can
be obtained by finding  elementary divisors of the matrix, that 
is, generators
of ideals generated by minors of $A_k$ of codimension $s$.
Elementary divisors can now be found by routine induction on $k$ 
to yield $D_s= gcd(W_{i,s})$, where $W_{i,s}$ 
is the set of products obtained from
$n_1...n_{i-1}m_in_{i+1}...n_k$ by dropping $s$ letters from it.
\end{proof}
The special form of the matrix representing  $H_{1}(M_L^{(2)}),$ where
$L=M(\frac{m_1}{n_1},...,\frac{m_k}{n_k})$, allows us to prove the GKH
 conjecture for such links.\\

{\bf Theorem 9}
{\it The GKH conjecture holds for
all alternating Montesinos links,
$L=M(\frac{m_1}{n_1},...,\frac{m_k}{n_k})$.}

\begin{proof}
We will show that no two different rational blocks (tangles) 
share the same element of the homology group. Different arcs inside 
each block represent different labels (see Fig.2.1 and the paragraph 
discussing it).
It is enough to compare the first block with the $j$th  block,
 where $1 < j <k$.
Assume now that an arc of the first block
represents the same homology as an arc of the $j$th block. Then,
in $H_1(M_ L^{(2)},\Z)$, we have for some $a$ and $b$, $0\leq a \leq m_1+n_1$
and
$0\leq b \leq n_j$:\\
$az_1=bz_j+x_{j-1}$ or equivalently
$az_1= bz_j+m_{j-1}z_{j-1} + m_{j-2}z_{j-2} +...+m_2z_2 + m_1z_1$,
and further
$$(a-m_1)z_1 -m_2z_2 -...- m_{j-2}z_{j-2} - m_{j-1}z_{j-1} -  bz_j=0.$$
Adding this relation
to the matrix of relations should keep the group unchanged. If the relation
$n_1z_1- n_jz_j$ is deleted, the group cannot be smaller. Therefore if we
replace the row $n_1z_1- n_jz_j$ by the row
$(a-m_1)z_1 -m_2z_2 -...- m_{j-2}z_{j-2} - m_{j-1}z_{j-1} -  bz_j$
we cannot decrease the (absolute value) of the determinant.
On the other hand, now the matrix has the form:

 \begin{displaymath}
\left[ \begin{array}{cccccccc}
n_1 & -n_2 & 0    & 0    &...&0   &... & 0   \\
n_1 & 0    & -n_3 & 0    &...&0   &... & 0   \\
n_1 & 0    & 0    & -n_4 &...&0   &... & 0   \\
... &...   &...   &...   &...&... &... &...  \\
a-m_1&-m_2 &-m_3  &-m_4  &...&-b  &... & 0   \\
... &...   &...   &...   &...&... &... &...  \\
n_1 & 0    & 0    & 0    &...&0   &... &-n_k \\
m_1 & m_2  & m_3  & m_4  &...&m_j &... &m_k
\end{array} \right]
\end{displaymath}
One can check that the absolute value of its determinant is smaller 
than $D$ of Proposition 8, unless $a=m_1+n_1,b=n_j$ and $j=2$. 
To demonstrate this, one uses 
properties of the matrix $A_k$ and its blocks of codimension $1$, and the 
fact that $0 < m_i < n_i$ (or $n_i=m_i=1$).
\end{proof}

\section{Future directions}

We expect that the method applied to prove the Generalized Kauffman-Harary 
Conjecture for Montesinos links can be extended to the case of 
$2$-algebraic links (i.e. algebraic links in the sense 
of Conway) and also to closed 3-braids. However, the general case 
requires new ideas. Exploiting the connection of the GKH conjecture 
to incompressible surfaces in the way outlined below seems to be a 
promising idea.

Assume that the GHK conjecture fails for an irreducible alternating link 
$L$ in $S^3.$  Then there are different arcs of its diagram labeled by 
$y_i$ and $y_j$ such 
that the element $y_iy_j^{-1}$ is homologically trivial in the double 
cover of $S^3$ branched along $L$. 
As the first step we analyze the possibility that this element is
homologically trivial in the unbranched double cover $\tilde M$ of
$S^3 -L$. In this case $y_iy_j^{-1}$ bounds a connected surface $\tilde F$
in $\tilde M$. Let $F$ be the projection of $\tilde F$ into $M.$
In order to show that $y_iy_j^{-1}$ is homologically non-trivial 
it is sufficient to prove that $F$ contains a meridional curve.
One may hope to show that by generalizing a theorem of Menasco
stating that every closed, incompressible surface in the 
exterior of an irreducible alternating knot contains a meridional curve,
\cite{M-1,M-2}.
This approach will be discussed in the sequel to this paper \cite{APS}.
\\
\ \\


\noindent
{\bf References} 
\begin{thebibliography}{99 1000-2000}

\bibitem{APS}
M.~M.~Asaeda, J.~H.~Przytycki, A.~S.~Sikora, in preparation.

\bibitem{Co}
J.~H.~Conway, An enumeration of knots and links,
{\nineit Computational problems in abstract algebra} (ed. J.~Leech),
Pergamon Press (1969) 329 - 358.

\bibitem{H-K}
F.~Harary, L.~H.~Kauffman, {\nineit Knots and Graphs I -- Arc Graphs 
and Colorings,} Advances in Applied Math. {\ninebf 22} (1999) 312--337.\\
\texttt{http://www.math.uic.edu/$\sim$kauffman/Papers.html}

\bibitem{Ka}
L.~H.~Kauffman, {\nineit Virtual Knot Theory,} talk at the special session
of AMS ``Algebraic Topology Based on Knots,"
AMS national meeting, January 2003, Baltimore, MD.

\bibitem{K-L}
L.~H.~Kauffman, S.~Lambropoulou,
{\nineit On the classification of rational tangles,} 
Advances in Applied Mathematics, to appear,
{http://users.ntua.gr/sofial/\#Preprints}

\bibitem{M-1}
W.~Menasco, {\nineit Closed incompressible surfaces in alternating knot 
and link complements,} Topology {\ninebf 23}(1) (1984) 37--44.
 
\bibitem{M-2}
W.~Menasco, {\nineit Determining incompressibility of surfaces in alternating 
knot and link complements,} Pacific J. Math. {\ninebf 117}(2) (1985) 353--370. 

\bibitem{PDDGS}
L.~Person, M.~Dunne, J.~DeNinno, B.~Guntel, L.~Smith,
{\nineit Colorings of rational, alternating knots and links,} preprint 2002.

\end{thebibliography}
\end{document}